\documentclass[12pt]{amsart}
\usepackage[T2A]{fontenc}
\usepackage[cp1251]{inputenc}
\usepackage[english]{babel}
\usepackage{amsmath,latexsym,amsthm,amsfonts,tipa,upgreek}
\usepackage{amssymb,amscd,graphpap, stmaryrd}

\newcommand\B{{\mathcal B}}

\theoremstyle{definition}

\begin{document}

\title{Counting coarse subsets of a countable group}
\author{Igor Protasov, Ksenia Protasova}
\date{}
\address{Department of Computer Science and Cybernetics, Kyiv University, Volodymyrska 64, 01033, Kyiv, Ukraine}
\email{i.v.protasov@gmail.com; }
\address{Department of Computer Science and Cybernetics, Kyiv University, Volodymyrska 64, 01033, Kyiv, Ukraine}
\email{ksuha@freenet.com.ua}
\maketitle

\begin{abstract}
For every countable group $G$, there are $2^{\omega}$ distinct classes of coarsely equivalent subsets of $G$.

\vspace{3 mm}

\vspace{3 mm}


\end{abstract}
\vspace{3 mm}

\section{Introduction and results}

Following \cite{b5}, \cite{b6}, we say that a {\em ball structure} is a triple $\B=(X,P,B)$, where $X$, $P$ are non-empty sets, and for all $x\in X$ and $\alpha\in P$, $B(x, \alpha)$ is a subset of $X$ which is called a {\em ball of radius} $\alpha$ around $x$. It is supposed that $x\in B(x, \alpha)$ for all $x\in X$, $\alpha\in P$.
The set $X$ is called the {\em support} of $\B$, $P$  is called the {\em set of radii}.

Given any $x\in X$, $A\subseteq X$, $\alpha\in P$, we set
$$B^*(x,\alpha)=\{y\in X:x\in B(y,\alpha)\},\ B(A,\alpha)=\bigcup_{a\in A}B(a,\alpha),\ B^*(A,\alpha)=\bigcup_{a\in A}B^*(a,\alpha).$$

A ball structure $\mathcal{B}=(X,P,B)$ is called a {\it ballean} if
\vspace{3 mm}

\begin{itemize}
\item{} for any $\alpha,\beta\in P$, there exist $\alpha',\beta'$ such that, for every $x\in X$,
$$B(x,\alpha)\subseteq B^*(x,\alpha'),\ B^*(x,\beta)\subseteq B(x,\beta');$$

\item{} for any $\alpha,\beta\in P$, there exists $\gamma\in P$ such that, for every $x\in X$,
$$B(B(x,\alpha),\beta)\subseteq B(x,\gamma);$$

\item{} for any $x,y\in X$, there exists $\alpha\in P$ such that $y\in B(x, \alpha)$.
\end{itemize}
\vspace{3 mm}

We note that a ballean can be considered as an asymptotic counterpart of a uniform space, and could be defined \cite{b7} in terms of  entourages of the diagonal $\Delta_X$ in $X\times X$. In this case a ballean is called a {\em coarse structure}. For categorical look at the balleans and coarse structures as "two faces of the same coin" see \cite{b2}.

Let $\mathcal{B}=(X,P,B)$, $\mathcal{B'}=(X',P',B')$ be balleans.
A mapping $f:X\to X'$ is called  {\it coarse} if, for every $\alpha\in P$, there exists $\alpha'\in P'$ such that, for every $x\in X$, $f(B(x,\alpha))\subseteq B'(f(x),\alpha')$.

A bijection $f:X\rightarrow X'$ is called an {\it asymorphism} between $\B$ and $\B'$ if $f$ and $f^{-1}$ are coarse. In this case $\B$ and $\B'$ are called {\it asymorphic}.

Let $\B = (X,P,B)$ be a ballean. Each subset $Y$ of $X$ defines a {\it subballean} $\B_Y = (Y,P,B_Y)$, where $B_Y(y,\alpha) = Y \cap B(y, \alpha)$. A subset $Y$ of $X$ is called {\it large} if $X = B(Y, \alpha)$, for $\alpha \in P$. Two balleans $\B$ and $\B'$ with supports $X$ and $X'$ are called {\it coarsely equivalent} if there exist large subsets $Y\subseteq X$ and $Y' \subseteq X'$ such that the subballeans $\B_Y$ and $\B'_{Y'}$ are asymorphic.

Every infinite group $G$ can be considered as the ballean $(G, \mathfrak{F}_{G}, B)$, where $ \mathfrak{F}_{G} $ is the family of all finite subsets of $G$,  $B(g,F)= Fg \bigcup \{g\}$.

We note that finitely generated groups are finitary coarsely equivalent if and only if $G$ and $H$ are quasi-isometric \cite[Chapter 4]{b3}.

A classification of countable locally finite groups (each finite subset generates finite subgroup) up to  asymorphisms is obtained in \cite{b4}(see also \cite[p. 103]{b5}).
\vspace{3 mm}

{\it Two countable locally finite groups $G_{1}$ and $G_{2}$ are  asymorphic if and only if the following conditions hold:
\vspace{3 mm}

$(i)$ for every finite subgroup $F\subset G_{1}$, there exists a finite subgroup $H$ of $G_{2}$ such that $|F|$ is a divisor of $|H|$;
\vspace{3 mm}

$(ii)$ for every finite subgroup $H$ of $G_{2}$, there exists a finite subgroup $F$ of $G_{1}$ such that $|F|$ is a divisor of $|F|$.}
\vspace{3 mm}

It follows that there are continuum many distinct types of countable locally finite groups and each group is asymorphic to some direct sum of finite cyclic groups.

The following coarse classification of countable Abelian groups is obtained in \cite{b1}.
\vspace{3 mm}

{\it Two countable Abelian groups are  coarsely equivalent if  and only if the torsion-free ranks of $G$  and $H$
coincide and  $G$  and $H$ are  either  both finitely generated or infinitely generated.}
\vspace{3 mm}

In particular, any two countable torsion Abelian groups are  coarsely equivalent.

Given a group $G$, we consider each non-empty subsets as a subballean of $G$ and say that a class  of
all pairwise coarsely equivalent subsets  is a {\it coarse subset } of $G$.

For a countable group $G$, we prove that there as many coarse subsets of $G$  as possible by the cardinal arithmetic.

{\bf Theorem.} {\it 
For a countable group $G$, there are $2^{\omega}$  coarse subsets of $G$.}

Every countable group $G$ contains either countable finitely generated subgroup or countable locally finite subgroup, so we split the proof into corresponding cases.

\section{Proof: finitely generated case}

2.1. We take a finite system $S,$  $S=S^{-1}$ of generators of $G$ and consider the Cayley graph $\Gamma$ with the set of vertices $G$ and the set of edges $\{\{g,h\}: gh^{-1} \in S,  \  g\neq h\}$.  We denote by $\rho$ the path metric on $\Gamma$  and choose a geodesic ray $V=\{v_{n}: n\in\omega\}$, $v_{0}$ is the identity of  $G$, $\rho(v_{n}, v_{m})=|n-m|$.

Then the subballean of $G$ with the support $V$ is asymorphic to the metric ballean $ \  (\mathbb{N}, \ \mathbb{N}\bigcup\{0\},  \  B)$,  where $B(x,r)=\{y\in\mathbb{N}: d(x,y )\leq r\},  \  d(x,y)= |x-y|$. Thus, it suffices to find a family  $\mathfrak{F}$, $|\mathfrak{F}|= 2^{\omega}$  of pairwise coarsely non-equivalent subsets of $\mathbb{N}$.
\vspace{3 mm}

2.2. We  choose a sequence $(I_{n})_{n\in\omega}$ of intervals of $\mathbb{N}$,  $I_{n}=[a_{n}, b_{n}]$,  $  \  b_{n}<a_{n+1}$ such that
\vspace{3 mm}

$(1)  \   \   \   \    b_{n}- a_{n} > n \  a_{n}.$
\vspace{3 mm}

Then we take an almost disjoint family $\mathcal{A}$ of infinite
subsets of $\omega$ such that $|\mathcal{A}|=2^{\omega}$. Recall that $\mathcal{A}$  is almost disjoint if $|W \bigcap W^{\prime}|<\omega$ for all distinct $W, \  W^{\prime} \in  \mathcal{A}$.

For each $W\in\mathcal{A}$, we denote $I_{W}= \bigcup\{I_{n}: n\in W\}$. To show that $\mathfrak{F}=\{ I_{W}: W\in \mathcal{A}\}$ is the desired family of subsets of $\mathbb{N}$, we take distinct $W, W^{\prime}\in\mathcal{A}$ and assume that $I_{W}$, $I_{W^{\prime}}$ are coarsely equivalent. Then there exist large subsets $X, X^{\prime}$  of $I_{W}$, $I_{W^{\prime}}$, and an asymorphism $f: X\longrightarrow X^{\prime}$. We choose $r\in\mathbb{N}$ such  that
$$I_{W}\subseteq B(X, r),  \   I_{W^{\prime}}\subseteq B(X^{\prime}, r) $$ and note that if an interval $I$  of length $2r$ is contained in $I_{W}$  then $I$  must contain at least one point of $X$, and the same holds for the pair
$I_{W^{\prime}}, X^{\prime}$.

Since $f$ is an asymorphism, we can  take $t\in\mathbb{N}$  such that, for all $x\in X$, $ \ x^{\prime}\in X^{\prime}$,
\vspace{3 mm}

$(2)  \   \   \   \   f(B_{X}(x, 2r+2)\subseteq  B_{X^{\prime}}(f(x), t);$
\vspace{3 mm}

$(3)  \   \   \   \    f^{-1}(B_{X^{\prime}}(x^{\prime}, 2r+2)\subseteq  B_{X}(f^{-1}(x^{\prime}), t).$
\vspace{3 mm}

We use $(1)$ to choose $m\in W\backslash W^{\prime}$, $m> \max (W\bigcap W^{\prime})$ such that

$(4)  \   \   \   \   b_{m}-a_{m} >  2r   a_{m};$
\vspace{3 mm}

$(5)  \   \   \   \    b_{m}-a_{m} >  2t.$
\vspace{5 mm}

2.3. We denote $Z=X\bigcap [a_{m}, b_{m}]$  and enumerate $Z$  in increasing order $Z=\{z_{0}, \ldots, z_{k}\}$.  Then $d(z_{i}, z_{i+1})\leq 2r + 2$  because otherwise the interval $[z_{i}+1, z_{i+1}-1]$  of length $2r$ has no points of
 $X$.

If $f(z_{0})< a_{m}$ then, by $(2)$ and $(5)$,  $f(Z)\subseteq [1, a_{m}-1]$.
On the other hand, $k\geq(b_{m}-a_{m})/2r - 1$ and, by $(4)$, $  \  (b_{m}-a_{m})/2r > a_{m}$.  Hence, $k >a_{m}-1$ contradicting $f(Z)\subseteq [1, a_{m}-1]$ because $f$ is a bijection.

If $f(z_{0})> b_{m}$ then  we take $s\in W^{\prime}$ such that $f(z_{0})\in I_{s}$. Since $m> \max (W\wedge  W^{\prime})$ and $s> m$,  we have  $s\in  W^{\prime}\setminus W$, so we can repeat above argument for $f^{-1}$ and $I_{s}$ in place of $f$  and $I_{m}$  with usage $(3)$ instead of $(2)$.

\section{ Proof: locally finite case}

3.1. Let $G$ be an arbitrary countable group and let $X$, $A$ be infinite subsets of $G$. Suppose that there exist an infinite subset $Y$  of $X$, a partition $A= B\bigcup C$  and $k, l\in\mathbb{N} $,  $k <  l$ such that
\vspace{3 mm}

$(6)  \   \   \   \ $  there exists $ H\in \mathfrak{F}_{G}$ such that, for every $y\in Y$,
$$|B_{X}(y, H)|\geq k;$$
\vspace{3 mm}

$(7)  \   \   \   \ $  for every  $ F  \in \mathfrak{F}_{G}$, there exists $ Y^{\prime}\in \mathfrak{F}_{G}$ such that, for every $y\in Y\backslash  Y^{\prime}$,
$$|B_{X}(y, H)|\geq  l;$$
\vspace{3 mm}

$(8)  \   \   \   \ $  there exists $ K\in \mathfrak{F}_{G}$ such that, for every $b\in B$,
$$|B_{A}(b, K)|>l;$$
\vspace{3 mm}

$(9)  \   \   \   \ $  for every  $ F  \in \mathfrak{F}_{G}$, there exists $ C^{\prime}\in \mathfrak{F}_{G}$ such that, for every $c\in C\backslash  C^{\prime}$,
$$|B_{A}(c, F)|<  k.$$
\vspace{3 mm}

Then $X$ and $A$ are not asymorphic.

We suppose the contrary and let $f: X\longrightarrow A$ be an asymorphism. We take an infinite subset $I$  of $Y$ such that either $f(I)\subset C$ or $f(I)\subset B$.

Assume that $f(I)\subset C$ and choose $ F  \in \mathfrak{F}_{G}$ such that, for every  $x\in X$,
$$f(B_{X} (x, H))\subseteq B_{A} (f(x), F). $$
For this $F$, we use $(9)$ to choose corresponding $C^{\prime}$.  We take $y\in I$ such that $f(y)\in C\setminus C^{\prime}$. By $(6)$, $f(B(y,H))\geq k$. By $(9)$, $ B_{A}(f(y),F)< k$  and we get a contradiction because $f$ is a bijection.

If $f(I)\subset B$  then, by $(8)$, $B_{A}(b,K)>l$ for every $b\in f(I)$. Since $f^{-1}$ is coarse, there is $ F  \in \mathfrak{F}_{G}$ such that, for every  $a\in A$
$$f^{-1}(B_{A}(a), K)\subseteq B_{X} (f^{-1}(a), F). $$
For this $F$, we choose $Y^{\prime}$ satisfying $(7)$  and get a contradiction.

\vspace{5 mm}

3.2.  Now we assume that $G$  is locally finite and show a plan  how to  choose the desired family $\mathfrak{F}$, $|\mathfrak{F}| = 2^{\omega}$ of pairwise coarsely non-equivalent subsets of $G$.

We construct some special sequence $(Y_{n})_{n\in\omega}$ of pairwise disjoint subsets of $G$. Then we take a family $\mathcal{A}$ of almost disjoint infinite subset of $\omega$, $|\mathfrak{F}| = 2^{\omega}$, denote
\vspace{3 mm}

$(10)  \   \   \   \  X_{W} =\bigcup\{ Y_{n}: n\in W\},  \  \  W\in \mathcal{A},$
\vspace{3 mm}

and get  $\mathfrak{F}$ as $\{ X_{W}: W\in \mathcal{A}\}.$

\vspace{5 mm}

3.3. We represent $G$ as the union of an increasing chain $\{F_{n}: n\in\omega\}$ of finite subgroups such that

$(11)  \   \   \   \  |F_{n+1}| >   |F_{n}|^{2}.$
\vspace{3 mm}

Then we choose a double sequence $(g_{nm})_{n,m\in\omega}$ of elements of $G$ such that
\vspace{3 mm}

$(12)  \   \   \   \  F_{n} F_{m} g_{n m} \bigcap F_{i} F_{j} g_{ij}=\emptyset$  for all distinct $(n, m)$, $(i, j)$ from $\omega\times\omega$,
and put $$Y_{n}=\bigcup\{ F_{m}g_{n m}:m\in\omega\}.$$

\vspace{5 mm}

3.4. We take distinct $W, W^{\prime}\in\mathcal{A}$ and prove that $X_{W}$  and $X_{W^{\prime}}$  (see (10)) are not coarsely equivalent. We suppose the contrary and choose large asymorphic subsets $Z_{W}$ and $Z_{W^{\prime}}$ of $X_{W}$  and  $X_{W^{\prime}}$. Then we take $t\in\omega$ such that  $$X_{W}\subseteq F_{t}  Z_{W},   \   X_{W^{\prime}}\subseteq F_{t}  Z_{W^{\prime}}. $$  If $n>t$
and either $F_{n}g_{nm} \subset X_{W}$  or $F_{n}g_{nm}\subset X_{W^{\prime}}$ then
\vspace{3 mm}

$(13)  \   \   \   \  |F_{n}  g_{n m} \bigcap Z_{W}|\geq \frac{|F_{n}|}{|F_{t}|},  \  \
 |F_{n} g_{nm}\bigcap Z_{W^{\prime}}| \geq \frac{|F_{n}|}{|F_{t}|}.$
 \vspace{3 mm}

 To apply 3.1, we choose $s\in W\setminus W^{\prime}, \  s>t $  and denote

 $$X= Z_{W},  \  \  Y=Y_{s}\bigcap Z_{W},  \  \  A=Z_{W^{\prime}},  \  \  B=\bigcup\{Y_{i}: i\in W^{\prime},  \  \  i>s\}, $$
$$C= \bigcup\{Y_{i}: i\in W^{\prime},  \  i<s \},  \  \  k=\frac{|F_{s}|}{|F_{t}|},  \  \  l=|F_{s}|. $$

By (13) with $s=n$,  we get (6). By (12) with $s=n$, we get (7).

If $n>s$  then $|F_{n}|/|F_{t}|> |F_{n}|/|F_{s}|$. By (11),   $|F_{n}|/|F_{s}|> |F_{s}|$, so $|F_{n}|/|F_{t}|> |F_{s}|$ and, by (13), we have (8).

If $n<s$ then $|F_{n}|<|F_{s}|/|F_{t}|$  and, by (12), we get (9).

 \section{ Comments}

 A subset $A$ of an infinite group $G$ is called

 \vspace{3 mm}

$\bullet$  {\it thick} if, for every $F\in\mathfrak{F}_{G}$ , there exists  $g\in A$ such that $Fg\subset A$;

\vspace{3 mm}

$\bullet$  {\it small} if $L\setminus A$ is large for every large subset $L$ of $G$;

\vspace{3 mm}

$\bullet$  {\it thin} if, for every $F\in\mathfrak{F}_{G}$ , there exists  $H\in \mathfrak{F}_{G}$ such that $B_{A} (g, F)= \{g\}$ for each $g\in A\setminus H$.

\vspace{3 mm}

A subset $A$ is thick if and only if $L\bigcap A\neq\emptyset$ for every large subset $L$ of $G$. For a countable group $G$, in the proof of Theorem, we construct $2^{\omega}$ pairwise coarsely non-equivalent thick subsets of $G$.

Every large subset $L$ of $G$ is coarsely equivalent to $G$, so any two large subsets of $G$ are coarsely equivalent. If $G$ is countable  then any two thin subset $S, T$ of $G$ are asymorphic: any bijection $f: S\longrightarrow T$ is an asymorphism. Every thin subset is small. But a small subset  $S$  of $G$ could be asymorphic to $G$: we take a group $G$ containing a subgroup $S$ isomorphic to $G$ such that the index of $S$ in $G$ is infinite.

\vspace{3 mm}
CONTACT INFORMATION

\end{document}